\documentclass[11pt,a4paper]{article}
\usepackage{amssymb}
\usepackage{amsfonts}
\usepackage{amsmath}
\usepackage{hyperref}
\usepackage{geometry}
\usepackage{relsize}

\input{tcilatex}
\begin{document}

\title{Error propagation in an explicit and an implicit numerical method for
Volterra integro-differential equations}
\author{J. S. C. Prentice \\
Senior Research Officer\\
Mathsophical Ltd.\\
Johannesburg, South Africa\\
Email: jpmsro@mathsophical.com}
\maketitle

\begin{abstract}
We study error propagation in both an explicit and an implicit method for
solving Volterra integro-differential equations. We determine the
relationship between local and global errors. We derive upper bounds for the
global error, and show that the global order for both methods is expected to
be first-order. A few numerical examples illustrate our results.
\end{abstract}

\section{Introduction}

Recently, we presented explicit and implicit numerical methods for solving
the Volterra integro-differential equation 
\begin{equation}
y^{\left( n\right) }\left( x\right) =f\left( x,y\right)
+\int\limits_{x_{0}}^{x}K\left( x,y\left( t\right) ,t\right) dt,\text{ \ \ \ 
}x>x_{0},  \label{problem}
\end{equation}%
using numerous examples to demonstrate the performance of the methods, and
also studying the stability of the methods \cite{prentice 2}\cite{prentice 1}%
. In this paper,\ we investigate the propagation of numerical error in these
methods. Not only is this an interesting study in its own right, it also
allows us to learn about upper bounds on the global error, and the order of
the error.

\section{Notation and terminology}

We deviate slightly form notation used in our previous work: here, $w$
denotes the approximate solution, and $y$ denotes the true solution. The
nodes are labelled as 
\begin{equation*}
x_{0}<x_{1}<x_{2}<\ldots <x_{i}<x_{i+1}<\ldots <x_{f}
\end{equation*}%
and $h$ is the uniform spacing between the nodes - the \textit{stepsize}. We
focus our attention on the case of $n=1$ in (\ref{problem}). We note that $f$
and $K$ are assumed to be suitably smooth so as to yield a unique solution
and, in particular, $K$ is not singular anywhere on the interval of
integration.

The explicit method is given by

\relsize{-1}%
\begin{align}
w_{i+1}=& \text{ }w_{i}+hf\left( x_{i},w_{i}\right)  \label{explicit} \\
& \text{ }+\frac{h^{2}}{2}\left( \sum\limits_{j=0}^{i}2K\left(
x_{i},w_{j},x_{j}\right) -K\left( x_{i},w_{0},x_{0}\right) -K\left(
x_{i},w_{i},x_{i}\right) \right)  \notag \\
& \equiv M_{E}\left( w_{i}\right) ,  \notag
\end{align}%
\relsize{1}where we have implicitly defined $M_{E}\left( w_{i}\right) .$

The implicit method is given by\relsize{-1}%
\begin{align}
w_{i+1}=\text{ }& w_{i}+hf\left( x_{i+1},w_{i+1}\right)  \label{implicit} \\
& \text{ }+\frac{h^{2}}{2}\left( \sum\limits_{j=0}^{i+1}2K\left(
x_{i+1},w_{j},x_{j}\right) -K\left( x_{i+1},w_{0},x_{0}\right) -K\left(
x_{i+1},w_{i+1},x_{i+1}\right) \right)  \notag \\
& \equiv M_{I}\left( w_{i}\right) ,  \notag
\end{align}%
\relsize{1}where we have implicitly defined $M_{I}\left( w_{i}\right) .$

We also define $M_{E}\left( y_{i}\right) $ and $M_{I}\left( y_{i}\right) $
as follows: \relsize{-1}%
\begin{align*}
y_{i+1}=& \text{ }y_{i}+hf\left( x_{i},y_{i}\right) \\
& \text{ }+\frac{h^{2}}{2}\left( \sum\limits_{j=0}^{i}2K\left(
x_{i},y_{j},x_{j}\right) -K\left( x_{i},y_{0},x_{0}\right) -K\left(
x_{i},y_{i},x_{i}\right) \right) \\
& \equiv M_{E}\left( y_{i}\right) .
\end{align*}%
\begin{align*}
y_{i+1}=\text{ }& y_{i}+hf\left( x_{i+1},y_{i+1}\right) \\
& +\frac{h^{2}}{2}\left( \sum\limits_{j=0}^{i+1}2K\left(
x_{i+1},y_{j},x_{j}\right) -K\left( x_{i+1},y_{0},x_{0}\right) -K\left(
x_{i+1},y_{i+1},x_{i+1}\right) \right) \\
& \equiv M_{I}\left( y_{i}\right) .
\end{align*}%
\relsize{1}

The \textit{global} error at $x_{i+1}$ is defined as 
\begin{align*}
\Delta _{i+1}& =w_{i+1}-y_{i+1}=M_{E}\left( w_{i}\right) -y_{i+1}\text{ \ \
(explicit method)} \\
\Delta _{i+1}& =w_{i+1}-y_{i+1}=M_{I}\left( w_{i}\right) -y_{i+1}\text{ \ \
(implicit method)}
\end{align*}%
and the \textit{local} error at $x_{i+1}$ is defined as%
\begin{align*}
\varepsilon _{i+1}& =M_{E}\left( y_{i}\right) -y_{i+1}\text{ \ \ (explicit
method)} \\
\varepsilon _{i+1}& =M_{I}\left( y_{i}\right) -y_{i+1}\text{ \ \ (implicit
method)}
\end{align*}%
The local error has the form%
\begin{equation*}
\varepsilon _{i+1}=\varepsilon _{i+1}^{D}+\varepsilon _{i+1}^{Q}=O\left(
h^{2}\right)
\end{equation*}%
where $\varepsilon _{i+1}^{D}$ is the error associated with the Euler
approximation to the derivative in the IDE, and $\varepsilon _{i+1}^{Q}$ is
the error associated with the composite Trapezium approximation to the
integral in the IDE. These are $O\left( h\right) ,$ at worst, but on
multiplication by $h,$ as required by the structure of the methods, these
errors acquire, at worst, an $O\left( h^{2}\right) $ character. The precise
form of these errors will not concern us here; it is enough for our purposes
to simply accept that $\varepsilon _{i+1}=O\left( h^{2}\right) .$
Nevertheless, we discuss this matter to some extent in the Appendix.

\section{Analysis - explicit case}

We consider the explicit case first, and do so in detail. The parameters $%
\xi $ and $\eta $ denote values appropriate for the various Taylor residual
terms that arise.

\subsection{Error propagation}

At $x_{1}$ we have\relsize{0}

\begin{align*}
w_{1}& =\text{ }w_{0}+hf\left( x_{0},w_{0}\right) \\
\Rightarrow \Delta _{1}+y_{1}& =\Delta _{0}+y_{0}+hf\left( x_{0},\Delta
_{0}+y_{0}\right) \\
& =\Delta _{0}+y_{0}+hf\left( x_{0},y_{0}\right) +\Delta _{0}hf_{y}\left(
x_{0},\xi _{0}\right) \\
\Rightarrow \Delta _{1}& =\left[ y_{0}+hf\left( x_{0},y_{0}\right) -y_{1}%
\right] +\Delta _{0}\left( 1+hf_{y}\left( x_{0},\xi _{0}\right) \right) \\
& =\left[ M_{E}\left( y_{0}\right) -y_{1}\right] +\Delta _{0}\left(
1+hf_{y}\left( x_{0},\xi _{0}\right) \right) \\
& =\varepsilon _{1}+\Delta _{0}\left( 1+hf_{y}\left( x_{0},\xi _{0}\right)
\right) .
\end{align*}%
\relsize{0}

At $x_{2}$ we have\relsize{-1}%
\begin{align*}
w_{2}=& \text{\thinspace }w_{1}+hf\left( x_{1},w_{1}\right) +\frac{h^{2}}{2}%
K\left( x_{0},w_{0},x_{0}\right) +\frac{h^{2}}{2}K\left(
x_{1},w_{1},x_{1}\right) \\
\Rightarrow \Delta _{2}+y_{2}=& \text{\thinspace }\Delta _{1}+y_{1}+hf\left(
x_{1},\Delta _{1}+y_{1}\right) +\frac{h^{2}}{2}K\left( x_{0},\Delta
_{0}+y_{0},x_{0}\right) +\frac{h^{2}}{2}K\left( x_{1},\Delta
_{1}+y_{1},x_{1}\right) \\
=& \,\Delta _{1}+y_{1}+hf\left( x_{1},y_{1}\right) +\Delta _{1}hf_{y}\left(
x_{1},\xi _{1}\right) \\
& +\frac{h^{2}}{2}\left( 
\begin{array}{c}
K\left( x_{0},y_{0},x_{0}\right) +\Delta _{0}K_{y}\left( x_{0},\eta
_{0},x_{0}\right) \ldots \\ 
\ldots +K\left( x_{1},y_{1},x_{1}\right) +\Delta _{1}K_{y}\left( x_{1},\eta
_{1},x_{1}\right)%
\end{array}%
\right) \\
\Rightarrow \Delta _{2}=& \left[ y_{1}+hf\left( x_{1},y_{1}\right) +\frac{%
h^{2}}{2}\left( K\left( x_{0},y_{0},x_{0}\right) +K\left(
x_{1},y_{1},x_{1}\right) \right) -y_{2}\right] \\
& +\Delta _{1}\left( 1+hf_{y}\left( x_{1},\xi _{1}\right) +\frac{h^{2}}{2}%
K_{y}\left( x_{1},\eta _{1},x_{1}\right) \right) +\Delta _{0}\frac{h^{2}}{2}%
K_{y}\left( x_{0},\eta _{0},x_{0}\right) \\
=& \text{\thinspace }\left[ M_{E}\left( y_{1}\right) -y_{2}\right] +\Delta
_{1}\left( 1+hf_{y}\left( x_{1},\xi _{1}\right) +\frac{h^{2}}{2}K_{y}\left(
x_{1},\eta _{1},x_{1}\right) \right) +\Delta _{0}\frac{h^{2}}{2}K_{y}\left(
x_{0},\eta _{0},x_{0}\right) \\
=& \text{\thinspace }\varepsilon _{2}+\Delta _{1}\left( 1+hf_{y}\left(
x_{1},\xi _{1}\right) +\frac{h^{2}}{2}K_{y}\left( x_{1},\eta
_{1},x_{1}\right) \right) +\sum\limits_{j=0}^{0}\Delta _{j}\frac{h^{2}}{2}%
K_{y}\left( x_{j},\eta _{j},x_{j}\right) ,
\end{align*}%
\relsize{1}and at $x_{3}$ we have\relsize{-1}%
\begin{align*}
w_{3}=& \,w_{2}+hf\left( x_{2},w_{2}\right) +\frac{h^{2}}{2}K\left(
x_{0},w_{0},x_{0}\right) +h^{2}K\left( x_{1},w_{1},x_{1}\right) +\frac{h^{2}%
}{2}K\left( x_{2},w_{2},x_{2}\right) \\
\Rightarrow \Delta _{3}+y_{3}=& \,\Delta _{2}+y_{2}+hf\left( x_{2},\Delta
_{2}+y_{2}\right) +\frac{h^{2}}{2}K\left( x_{0},\Delta
_{0}+y_{0},x_{0}\right) +h^{2}K\left( x_{1},\Delta _{1}+y_{1},x_{1}\right) \\
& +\frac{h^{2}}{2}K\left( x_{2},\Delta _{2}+y_{2},x_{2}\right) \\
=& \text{\thinspace }\Delta _{2}+y_{2}+hf\left( x_{2},y_{2}\right) +\Delta
_{2}hf_{y}\left( x_{2},\xi _{2}\right) \\
& +\frac{h^{2}}{2}\left( 
\begin{array}{c}
K\left( x_{0},y_{0},x_{0}\right) +\Delta _{0}K_{y}\left( x_{0},\eta
_{0},x_{0}\right) +2K\left( x_{1},y_{1},x_{1}\right) \ldots \\ 
\ldots +2\Delta _{1}K_{y}\left( x_{1},\eta _{1},x_{1}\right) +K\left(
x_{2},y_{2},x_{2}\right) +\Delta _{2}K_{y}\left( x_{2},\eta _{2},x_{2}\right)%
\end{array}%
\right)
\end{align*}%
\begin{align*}
\Rightarrow \Delta _{3}=& \left[ y_{2}+hf\left( x_{2},y_{2}\right) +\frac{%
h^{2}}{2}\left( K\left( x_{0},y_{0},x_{0}\right) +2K\left(
x_{1},y_{1},x_{1}\right) +K\left( x_{2},y_{2},x_{2}\right) \right) -y_{3}%
\right] \\
& +\Delta _{2}\left( 1+hf_{y}\left( x_{2},\xi _{2}\right) +\frac{h^{2}}{2}%
K_{y}\left( x_{2},\eta _{2},x_{2}\right) \right) +\Delta _{1}h^{2}K\left(
x_{1},\eta _{1},x_{1}\right) +\Delta _{0}\frac{h^{2}}{2}K_{y}\left(
x_{0},\eta _{0},x_{0}\right) \\
=& \text{\thinspace }\left[ M_{E}\left( y_{2}\right) -y_{3}\right] +\Delta
_{2}\left( 1+hf_{y}\left( x_{2},\xi _{2}\right) +\frac{h^{2}}{2}K_{y}\left(
x_{2},\eta _{2},x_{2}\right) \right) +\Delta _{1}h^{2}K\left( x_{1},\eta
_{1},x_{1}\right) \\
& +\Delta _{0}\frac{h^{2}}{2}K_{y}\left( x_{0},\eta _{0},x_{0}\right) \\
=& \text{\thinspace }\varepsilon _{3}+\Delta _{2}\left( 1+hf_{y}\left(
x_{2},\xi _{2}\right) +\frac{h^{2}}{2}K_{y}\left( x_{2},\eta
_{2},x_{2}\right) \right) +\sum\limits_{j=0}^{1}\Delta _{j}\frac{h^{2}}{2}%
K_{y}\left( x_{j},\eta _{j},x_{j}\right) \\
& +\sum\limits_{j=1}^{1}\Delta _{j}\frac{h^{2}}{2}K_{y}\left( x_{j},\eta
_{j},x_{j}\right) \\
=& \,\varepsilon _{3}+\Delta _{2}\left( 1+hf_{y}\left( x_{2},\xi _{2}\right)
+\frac{h^{2}}{2}K_{y}\left( x_{2},\eta _{2},x_{2}\right) \right)
+\sum\limits_{j=0}^{1}\Delta _{j}\frac{h^{2}}{2}K_{y}\left( x_{j},\eta
_{j},x_{j}\right) \\
& +\sum\limits_{j=1}^{1}\Delta _{j}\frac{h^{2}}{2}K_{y}\left( x_{j},\eta
_{j},x_{j}\right) .
\end{align*}%
\relsize{1}In general, for $i>1,$ we have%
\begin{align}
\Delta _{i+1}=& \text{\thinspace }\varepsilon _{i+1}+\Delta _{i}\left(
1+hf_{y}\left( x_{i},\xi _{i}\right) +\frac{h^{2}}{2}K_{y}\left( x_{i},\eta
_{i},x_{i}\right) \right)  \notag \\
& +\frac{h^{2}}{2}\left( \sum\limits_{j=0}^{i-1}\Delta _{j}K_{y}\left(
x_{j},\eta _{j},x_{j}\right) +\sum\limits_{j=1}^{i-1}\Delta _{j}K_{y}\left(
x_{j},\eta _{j},x_{j}\right) \right)  \notag \\
\Rightarrow \Delta _{i+1}=\text{\thinspace }& \widetilde{\varepsilon }_{i+1}+%
\widetilde{\alpha }_{i}\Delta _{i},  \label{delta = eps + alpha*delta}
\end{align}%
where%
\begin{align*}
\widetilde{\varepsilon }_{i+1}& \equiv \varepsilon _{i+1}+\frac{h^{2}}{2}%
\left( \sum\limits_{j=0}^{i-1}\Delta _{j}K_{y}\left( x_{j},\eta
_{j},x_{j}\right) +\sum\limits_{j=1}^{i-1}\Delta _{j}K_{y}\left( x_{j},\eta
_{j},x_{j}\right) \right) \\
& =\varepsilon _{i+1}+\frac{h^{2}}{2}\left( \sum\limits_{j=0}^{i-1}2\Delta
_{j}K_{y}\left( x_{j},\eta _{j},x_{j}\right) -\frac{\Delta _{0}}{2}%
K_{y}\left( x_{0},\eta _{0},x_{0}\right) \right) \\
& =\varepsilon _{i+1}+h^{2}\sum\limits_{j=1}^{i-1}\Delta _{j}K_{y}\left(
x_{j},\eta _{j},x_{j}\right) \text{ \ if }\Delta _{0}=0.
\end{align*}%
and%
\begin{equation*}
\widetilde{\alpha }_{i}\equiv 1+hf_{y}\left( x_{i},\xi _{i}\right) +\frac{%
h^{2}}{2}K_{y}\left( x_{i},\eta _{i},x_{i}\right) =1+h\left( f_{y}\left(
x_{i},\xi _{i}\right) +\frac{h}{2}K_{y}\left( x_{i},\eta _{i},x_{i}\right)
\right) .
\end{equation*}%
Note that $\widetilde{\varepsilon }_{i+1}$ is not a local error, but it is
convenient to combine the terms in this way, as we shall soon see. Also, we
can write%
\begin{align*}
\widetilde{\varepsilon }_{i+1}& =\left(
C_{i+1}^{1}+\sum\limits_{j=0}^{i-1}\Delta _{j}K_{y}\left( x_{j},\eta
_{j},x_{j}\right) \right) h^{2} \\
& \equiv \left( C_{i+1}^{1}+C_{i+1}^{2}\right) h^{2} \\
& =\widetilde{C}_{i+1}h^{2}
\end{align*}%
wherein the coefficients $C_{i+1}^{1},C_{i+1}^{2}$ and $\widetilde{C}_{i+1}$
have been implicitly defined. Equation\ (\ref{delta = eps + alpha*delta}) is
the defining expression for the propagation of error in the explicit method.

In the remainder of this paper, we will assume $\Delta _{0}=0.$

\subsection{Upper bounds}

Assume $f_{y}+\frac{h}{2}K_{y}>0.$ With

\begin{align*}
\widetilde{\varepsilon }_{\max }& \equiv \max_{\left[ x_{0},x_{f}\right]
}\left\vert \widetilde{\varepsilon }_{i}\right\vert =\max_{\left[ x_{0},x_{f}%
\right] }\left\vert \widetilde{C}_{i}\right\vert h^{2}\equiv \widetilde{C}%
h^{2} \\
\widetilde{\alpha }& \equiv 1+\max_{\left[ x_{0},x_{f}\right] }\left( hf_{y}+%
\frac{h^{2}}{2}K_{y}\right) \equiv 1+hL \\
\Rightarrow L& =\max_{\left[ x_{0},x_{f}\right] }\left( f_{y}+\frac{h}{2}%
K_{y}\right) 
\end{align*}%
we find%
\begin{align}
\left\vert \Delta _{i+1}\right\vert & \leqslant \widetilde{\varepsilon }%
_{\max }\left( 1+\widetilde{\alpha }+\widetilde{\alpha }^{2}+\ldots +%
\widetilde{\alpha }^{i}\right)   \notag \\
& =\widetilde{\varepsilon }_{\max }\left( \frac{\widetilde{\alpha }^{i+1}-1}{%
\widetilde{\alpha }-1}\right)   \notag \\
& =\frac{\widetilde{\varepsilon }_{\max }}{hL}\left( \left( 1+hL\right)
^{i+1}-1\right)   \notag \\
& =\frac{\widetilde{\varepsilon }_{\max }}{hL}\left( \left( 1+\frac{\left(
i+1\right) hL}{i+1}\right) ^{i+1}-1\right)   \notag \\
& =\frac{\widetilde{C}h^{2}}{hL}\left( \left( 1+\frac{\left(
x_{i+1}-x_{0}\right) L}{i+1}\right) ^{i+1}-1\right)   \notag \\
& \approx \frac{\widetilde{C}h}{L}\left( e^{\left( x_{i+1}-x_{0}\right)
L}-1\right) \text{ \ for large }i.  \label{Ch/L}
\end{align}

If $f_{y}+\frac{h}{2}K_{y}<0$ we define $L\equiv -\max_{\left[ x_{0},x_{f}%
\right] }\left\vert f_{y}+\frac{h}{2}K_{y}\right\vert .$ We can then choose $%
h$ so that $\widetilde{\alpha }=1+hL>0,$ and we then find%
\begin{equation}
\left\vert \Delta _{i+1}\right\vert \lesssim \left\vert \frac{\widetilde{C}h%
}{L}\left( e^{\left( x_{i+1}-x_{0}\right) L}-1\right) \right\vert \approx
\left\vert -\frac{\widetilde{C}}{L}\right\vert h\text{ \ if }f_{y}+\frac{h}{2%
}K_{y}\ll 0.  \label{-C/hL}
\end{equation}

If $f_{y}+\frac{h}{2}K_{y}=0$ we have $\widetilde{\alpha }=1$ and so%
\begin{align*}
\left\vert \Delta _{i+1}\right\vert & \leqslant \widetilde{\varepsilon }%
_{\max }\underset{i+1\text{ times}}{\underbrace{\left( 1+1+1+\ldots
+1\right) }} \\
& =\widetilde{C}h\left( i+1\right) h \\
& =\widetilde{C}\left( x_{i+1}-x_{0}\right) h.
\end{align*}

We see that all of these bounds exhibit a first-order $\left( O\left(
h\right) \right) $\medskip\ character.

\subsection{Order}

\medskip Assume $h$ is sufficiently small so that $\widetilde{\alpha }%
_{i}\approx 1$ 
\begin{align*}
\Delta _{1}& =\widetilde{\varepsilon }_{1}+\widetilde{\alpha }_{0}\Delta
_{0}=\widetilde{\varepsilon }_{1}=\left( C_{1}^{1}+\Delta _{0}K_{y}\left(
x_{0},\eta _{0},x_{0}\right) \right) h^{2}=C_{1}^{1}h^{2} \\
\Delta _{2}& =\widetilde{\varepsilon }_{2}+\widetilde{\alpha }_{1}\Delta
_{1}\approx \left( C_{2}^{1}+\Delta _{0}K_{y}\left( x_{0},\eta
_{0},x_{0}\right) +\Delta _{1}K_{y}\left( x_{1},\eta _{1},x_{1}\right)
\right) h^{2}+\Delta _{1} \\
& =\left( C_{2}^{1}+C_{1}^{1}h^{2}K_{y}\left( x_{1},\eta _{1},x_{1}\right)
\right) h^{2}+C_{1}^{1}h^{2}\approx \left( C_{2}^{1}+C_{1}^{1}\right) h^{2}
\\
& =\left( \frac{C_{2}^{1}+C_{1}^{1}}{2}\right) 2h^{2}=\left( \frac{%
C_{2}^{1}+C_{1}^{1}}{2}\right) \left( 2h\right) h \\
\Delta _{3}& =\left( \frac{C_{3}^{1}+C_{2}^{1}+C_{1}^{1}}{3}\right)
3h^{2}=\left( \frac{C_{3}^{1}+C_{2}^{1}+C_{1}^{1}}{3}\right) \left(
3h\right) h \\
& \vdots \\
\Delta _{i+1}& =\left( \frac{C_{i+1}^{1}+\ldots +C_{2}^{1}+C_{1}^{1}}{i+1}%
\right) \left( i+1\right) h^{2}=\left( \frac{C_{i+1}^{1}+\ldots
+C_{2}^{1}+C_{1}^{1}}{i+1}\right) \left( \left( i+1\right) h\right) h
\end{align*}

\medskip Now, let $x_{d}>x_{0}$ and choose $h_{1}$ so that%
\begin{equation*}
x_{d}=x_{0}+m_{1}h_{1}.
\end{equation*}%
Hence, 
\begin{align*}
\Delta _{m_{1}}& =\left( \frac{C_{m_{1}}^{1}+\ldots +C_{2}^{1}+C_{1}^{1}}{%
m_{1}}\right) \left( m_{1}h_{1}\right) h_{1} \\
& =\left( \frac{C_{m_{1}}^{1}+\ldots +C_{2}^{1}+C_{1}^{1}}{m_{1}}\right)
\left( x_{d}-x_{0}\right) h_{1}.
\end{align*}%
Now, choose $m_{2}\neq m_{1}$ and $h_{2}$ such that%
\begin{equation*}
x_{d}=x_{0}+m_{2}h_{2}.
\end{equation*}%
Note that%
\begin{equation*}
m_{1}h_{1}=m_{2}h_{2}\Rightarrow h_{2}=\left( \frac{m_{1}}{m_{2}}\right)
h_{1}.
\end{equation*}%
Hence, 
\begin{align*}
\Delta _{m_{2}}& =\left( \frac{C_{m_{2}}^{1}+\ldots +C_{2}^{1}+C_{1}^{1}}{%
m_{2}}\right) \left( m_{2}h_{2}\right) h_{2} \\
& =\left( \frac{C_{m_{2}}^{1}+\ldots +C_{2}^{1}+C_{1}^{1}}{m_{2}}\right)
\left( x_{d}-x_{0}\right) h_{2} \\
& =\left( \frac{C_{m_{2}}^{1}+\ldots +C_{2}^{1}+C_{1}^{1}}{m_{2}}\right)
\left( x_{d}-x_{0}\right) \left( \frac{m_{1}}{m_{2}}\right) h_{1} \\
& \approx \Delta _{m_{1}}\left( \frac{m_{1}}{m_{2}}\right)
\end{align*}%
so that the global error at $x_{d}$ scales in the same way as the stepsize,
i.e. the global error is first-order. This aligns with the $O\left( h\right) 
$ nature of the upper bounds considered earlier. Note that this implies that
the explicit method is \textit{convergent} ($\Delta \rightarrow 0$ as $%
h\rightarrow 0).$

\section{Analysis - implicit case}

\subsection{Error propagation}

For the implicit method we have%
\begin{equation*}
w_{1}=\text{ }w_{0}+hf\left( x_{1},w_{1}\right) +\frac{h^{2}}{2}\left(
K\left( x_{1},w_{0},x_{0}\right) +K\left( x_{1},w_{1},x_{1}\right) \right)
\end{equation*}%
which gives (with $\Delta _{0}=0)$%
\begin{align*}
\Delta _{1}+y_{1}& =\text{ }y_{0}+hf\left( x_{1},\Delta _{1}+y_{1}\right) +%
\frac{h^{2}}{2}\left( K\left( x_{1},y_{0},x_{0}\right) +K\left( x_{1},\Delta
_{1}+y_{1},x_{1}\right) \right) \\
\Rightarrow \Delta _{1}& =\frac{\left[ y_{0}+hf\left( x_{1},y_{1}\right) +%
\frac{h^{2}}{2}\left( K\left( x_{1},y_{0},x_{0}\right) +K\left(
x_{1},y_{1},x_{1}\right) \right) -y_{1}\right] }{1-hf_{y}\left( x_{1},\xi
_{1}\right) -\frac{h^{2}}{2}K_{y}\left( x_{1},\eta _{1},x_{1}\right) } \\
\Rightarrow \Delta _{1}& =\widetilde{\widetilde{\alpha }}_{1}\left(
M_{I}\left( y_{0}\right) -y_{1}\right) \\
& =\widetilde{\widetilde{\alpha }}_{1}\varepsilon _{1},
\end{align*}%
wherein we have implicitly defined $\widetilde{\widetilde{\alpha }}_{1}.$

At $x_{2}$ we find%
\begin{equation*}
\Delta _{2}=\frac{\varepsilon _{2}+\Delta _{1}\left( 1+h^{2}K_{y}\left(
x_{1},\eta _{1},x_{1}\right) \right) }{1-hf_{y}\left( x_{2},\xi _{2}\right) -%
\frac{h^{2}}{2}K_{y}\left( x_{2},\eta _{2},x_{2}\right) }
\end{equation*}%
and at $x_{3}$ we find%
\begin{equation*}
\Delta _{3}=\frac{\varepsilon _{3}+\Delta _{2}\left( 1+h^{2}K_{y}\left(
x_{2},\eta _{2},x_{2}\right) \right) +\Delta _{1}h^{2}K_{y}\left( x_{1},\eta
_{1},x_{1}\right) }{1-hf_{y}\left( x_{3},\xi _{3}\right) -\frac{h^{2}}{2}%
K_{y}\left( x_{3},\eta _{3},x_{3}\right) }.
\end{equation*}%
In general (for $i>1),$ we have%
\begin{align}
\Delta _{i+1}& =\frac{\varepsilon _{i+1}+\sum\limits_{j=1}^{i-1}\Delta
_{j}h^{2}K_{y}\left( x_{j},\eta _{j},x_{j}\right) +\Delta _{i}\left(
1+h^{2}K_{y}\left( x_{i},\eta _{i},x_{i}\right) \right) }{1-hf_{y}\left(
x_{i+1},\xi _{i+1}\right) -\frac{h^{2}}{2}K_{y}\left( x_{i+1},\eta
_{i+1},x_{i+1}\right) }  \notag \\
\Rightarrow \Delta _{i+1}& =\widetilde{\widetilde{\varepsilon }}_{i+1}+%
\widetilde{\widetilde{\alpha }}_{i}\Delta _{i},  \label{error prop implicit}
\end{align}%
where%
\begin{align*}
\widetilde{\widetilde{\varepsilon }}_{i+1}& \equiv \frac{\varepsilon
_{i+1}+\sum\limits_{j=1}^{i-1}\Delta _{j}h^{2}K_{y}\left( x_{j},\eta
_{j},x_{j}\right) }{1-hf_{y}\left( x_{i+1},\xi _{i+1}\right) -\frac{h^{2}}{2}%
K_{y}\left( x_{i+1},\eta _{i+1},x_{i+1}\right) } \\
\widetilde{\widetilde{\alpha }}_{i}& \equiv \frac{1+h^{2}K_{y}\left(
x_{i},\eta _{i},x_{i}\right) }{1-hf_{y}\left( x_{i+1},\xi _{i+1}\right) -%
\frac{h^{2}}{2}K_{y}\left( x_{i+1},\eta _{i+1},x_{i+1}\right) }.
\end{align*}

\subsection{Upper bound}

We are most likely to use the implicit method when the IDE is stiff (both $%
f_{y}<0$ and $K_{y}<0).$ Hence, it is instructive to apply (\ref{error prop
implicit}) to the test equation \cite{prentice 1}%
\begin{align}
y^{\prime }\left( x\right) & =\lambda \left( y\left( x\right) -1\right)
+\gamma \int\limits_{0}^{x}y\left( t\right) dt  \label{test eqn} \\
y\left( 0\right) & =2,\text{ }\lambda <0,\text{ }\gamma <0  \notag
\end{align}%
(where $f=\lambda \left( y\left( x\right) -1\right) $ and $K=\gamma y\left(
t\right) ),$ with solution%
\begin{align*}
y\left( x\right) & =e^{m_{1}x}+e^{m_{2}x} \\
m_{1}& =\frac{\lambda -\sqrt{\lambda ^{2}+4\gamma }}{2},\text{ \ }m_{2}=%
\frac{\lambda +\sqrt{\lambda ^{2}+4\gamma }}{2}
\end{align*}%
when $m_{1}$ and $m_{2}$ are real $\left( \lambda ^{2}+4\gamma \geqslant
0\right) $, and%
\begin{equation*}
y\left( x\right) =2e^{\frac{\lambda x}{2}}\cos \left( \frac{\sqrt{\left\vert
\lambda ^{2}+4\gamma \right\vert }}{2}x\right)
\end{equation*}%
when $m_{1}$ and $m_{2}$ are complex $\left( \lambda ^{2}+4\gamma <0\right) $%
.

With $Z\equiv h\lambda $ and $W\equiv h^{2}\gamma ,$ we define $L$ by%
\begin{align*}
1+hL& \equiv \frac{1+h^{2}K_{y}}{1-hf_{y}-\frac{h^{2}}{2}K_{y}}=\frac{1+W}{%
1-Z-\frac{W}{2}} \\
\Rightarrow L& =\frac{Z+\frac{3W}{2}}{h\left( 1-Z-\frac{W}{2}\right) }.
\end{align*}%
Since $Z$ and $W$ are both negative, $L$ is negative, too. Also%
\begin{align*}
L=\frac{Z+\frac{3W}{2}}{h\left( 1-Z-\frac{W}{2}\right) }& =\frac{h\lambda +%
\frac{3h^{2}\gamma }{2}}{h\left( 1-h\lambda -\frac{h^{2}\gamma }{2}\right) }
\\
& =\frac{\lambda +\frac{3h\gamma }{2}}{\left( 1-h\lambda -\frac{h^{2}\gamma 
}{2}\right) }.
\end{align*}%
With

\begin{align*}
\widetilde{\widetilde{\varepsilon }}_{\max }& \equiv \max_{\left[ x_{0},x_{f}%
\right] }\left\vert \widetilde{\widetilde{\varepsilon }}_{i}\right\vert
\equiv \max_{\left[ x_{0},x_{f}\right] }\left\vert \widetilde{\widetilde{C}}%
_{i}\right\vert h^{2}\equiv \widetilde{\widetilde{C}}h^{2} \\
\widetilde{\widetilde{\alpha }}& \equiv 1+hL
\end{align*}%
we find%
\begin{align*}
\left\vert \Delta _{i+1}\right\vert & \leqslant \widetilde{\widetilde{%
\varepsilon }}_{\max }\left\vert 1+\widetilde{\widetilde{\alpha }}+%
\widetilde{\widetilde{\alpha }}^{2}+\ldots +\widetilde{\widetilde{\alpha }}%
^{i}\right\vert \\
& =\widetilde{\widetilde{\varepsilon }}_{\max }\left\vert \frac{\widetilde{%
\widetilde{\alpha }}^{i+1}-1}{\widetilde{\widetilde{\alpha }}-1}\right\vert
\\
& =\frac{\widetilde{\widetilde{\varepsilon }}_{\max }}{h\left\vert
L\right\vert }\left\vert \left( 1+hL\right) ^{i+1}-1\right\vert \\
& =\frac{\widetilde{\widetilde{\varepsilon }}_{\max }}{h\left\vert
L\right\vert }\left\vert \left( 1+\frac{\left( i+1\right) hL}{i+1}\right)
^{i+1}-1\right\vert \\
& =\frac{\widetilde{\widetilde{C}}h^{2}}{h\left\vert L\right\vert }%
\left\vert \left( 1+\frac{\left( x_{i+1}-x_{0}\right) L}{i+1}\right)
^{i+1}-1\right\vert \\
& \approx \left\vert \frac{\widetilde{\widetilde{C}}h}{L}\left( e^{\left(
x_{i+1}-x_{0}\right) L}-1\right) \right\vert \text{ \ for large }i.
\end{align*}%
Since $L<0,$ we note that%
\begin{equation*}
\left\vert \frac{\widetilde{\widetilde{C}}h}{L}\left( e^{\left(
x_{i+1}-x_{0}\right) L}-1\right) \right\vert \rightarrow \left\vert -\frac{%
\widetilde{\widetilde{C}}h}{L}\right\vert
\end{equation*}%
if $L\ll 0$, and/or if $x_{i+1}-x_{0}$ becomes large (similar to the case
considered in (\ref{-C/hL})).

\subsection{Order}

To analyze the order of the implicit method, we assume that $h$ is small
enough so that%
\begin{align*}
\widetilde{\widetilde{\varepsilon }}_{i+1}& \approx \varepsilon
_{i+1}+\sum\limits_{j=1}^{i-1}\Delta _{j}h^{2}K_{y}\left( x_{j},\eta
_{j},x_{j}\right) \\
\widetilde{\widetilde{\alpha }}_{i}& \approx 1.
\end{align*}%
Similar reasoning to the explicit case can now be used to find that the
implicit method is expected to be first-order. Furthermore, this implies
that the implicit method is convergent.

\section{Comments}

For the explicit method, given that $\Delta _{1}=\varepsilon
_{1}=\varepsilon _{1}^{D}+\varepsilon _{1}^{Q},$ and given that all
subsequent global errors are written in terms of local errors and prior
global errors, we have that $\Delta _{i}$ is a function of local errors $%
\varepsilon ^{D}$ and $\varepsilon ^{Q},$ Jacobians $f_{y}$ and $K_{y}$ and
the stepsize $h$. For example, we find%
\begin{equation*}
\Delta _{4}=\varepsilon _{4}+\widetilde{\alpha }_{3}\varepsilon _{3}+%
\widetilde{\alpha }_{3}\widetilde{\alpha }_{2}\varepsilon _{2}+\widetilde{%
\alpha }_{3}\widetilde{\alpha }_{2}\widetilde{\alpha }_{1}\varepsilon
_{1}+h^{2}\left( \varepsilon _{1}K_{y}^{1}+\left( \varepsilon _{2}+%
\widetilde{\alpha }_{1}\varepsilon _{1}\right) K_{y}^{2}\right)
\end{equation*}%
where%
\begin{align*}
\widetilde{\alpha }_{k}& =1+hf_{y}\left( x_{k},\xi _{k}\right) +\frac{h^{2}}{%
2}K_{y}\left( x_{k},\eta _{k},x_{k}\right) \\
K_{y}^{k}& \equiv K_{y}\left( x_{k},\eta _{k},x_{k}\right) .
\end{align*}%
Similar expressions obtain for $\Delta _{5,}\Delta _{6}$ and so on, and also
for the case of the implicit method. It is interesting to note that, if the
global $\Delta _{i}$ error is known and the Jacobians can be reliably
estimated (such as for the test equation), then the local errors $%
\varepsilon _{i}$ (for the explicit method) can be estimated via the sequence%
\begin{align}
\varepsilon _{1}& =\Delta _{1},  \notag \\
\varepsilon _{2}& =\Delta _{2}-\widetilde{\alpha }_{1}\Delta _{1},
\label{sequence} \\
\varepsilon _{3}& =\Delta _{3}-\widetilde{\alpha }_{2}\Delta
_{2}-h^{2}\Delta _{1}K_{y}^{1},  \notag \\
\varepsilon _{4}& =\Delta _{4}-\widetilde{\alpha }_{3}\Delta
_{3}-h^{2}\Delta _{1}K_{y}^{1}-h^{2}\Delta _{2}K_{y}^{2}  \notag
\end{align}%
and so on. A similar sequence can be found for the implicit method.

\section{Numerical examples}

A few simple examples, using the test equation, will serve to illustrate
some of the aspects of our analysis.

\begin{enumerate}
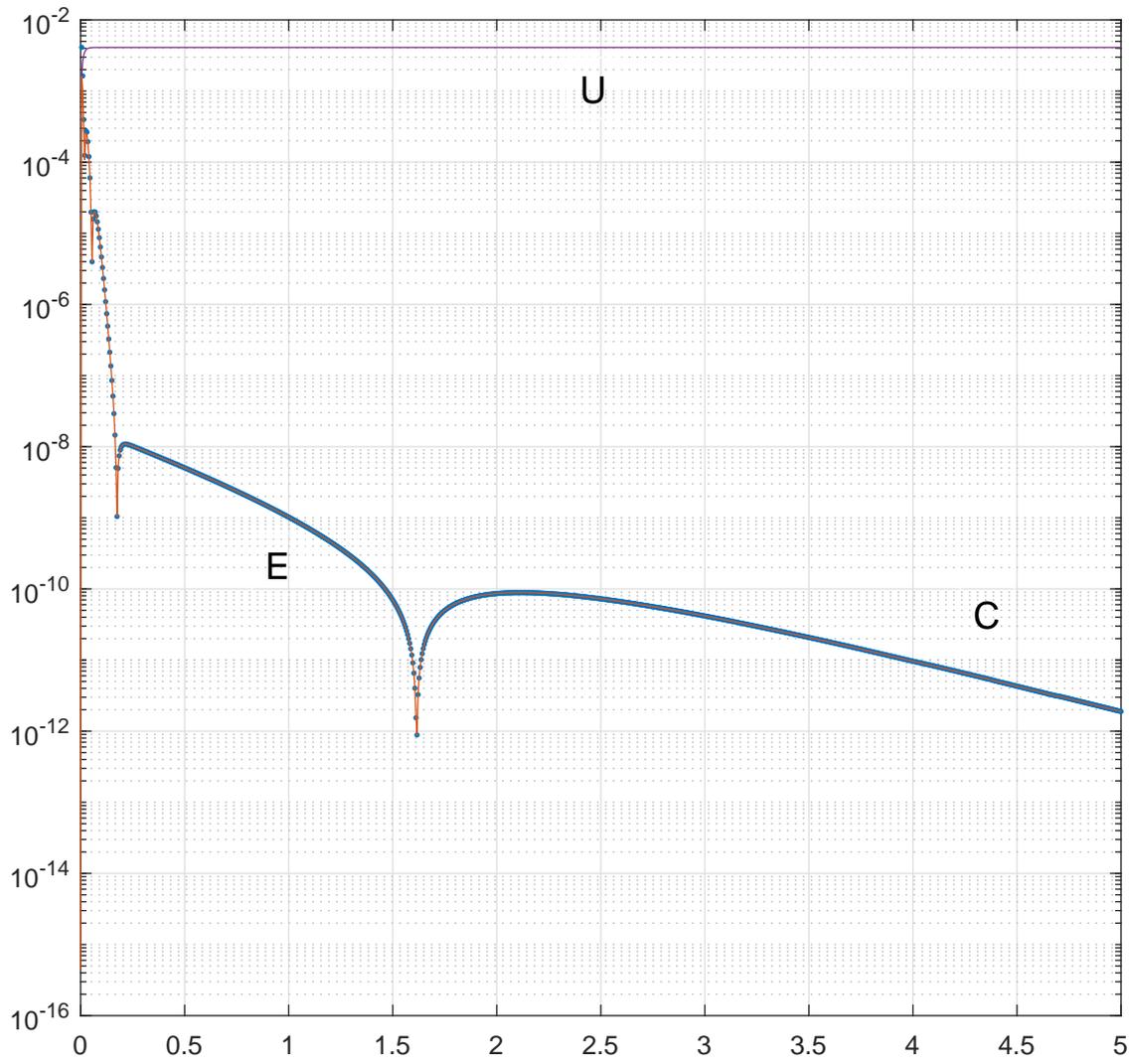
\item \textbf{Figure 1}. Here, we solve (\ref{test eqn}) with $\lambda
=-100,\gamma =-200$ and $h=5\times 10^{-3}$ using the explicit method. The
stepsize is small enough to ensure a stable solution. We show $\left\vert
\Delta _{i}\right\vert $ (the solid red line, labelled E), and the quantity $%
\left\vert \frac{\widetilde{C_{i}}h}{L}\right\vert $ determined using (\ref%
{Ch/L}), i.e.%
\begin{equation*}
\left\vert \frac{\widetilde{C}_{i}h}{L}\right\vert =\frac{\left\vert \Delta
_{i}\right\vert }{\left\vert e^{\left( x_{i+1}-x_{0}\right) L}-1\right\vert }%
.
\end{equation*}%
We indicate $\left\vert \frac{\widetilde{C_{i}}h}{L}\right\vert $ with the
blue dots (labelled C) which appear to be superimposed on the curve for $%
\left\vert \Delta _{i}\right\vert .$ This is due to the fact that $L=f_{y}+%
\frac{h}{2}K_{y}=\lambda +\frac{h}{2}\gamma \ll 0,$ and so $\left\vert
e^{\left( x_{i+1}-x_{0}\right) L}-1\right\vert \approx 1.$ From this curve
we estimate $\max \left\vert \frac{\widetilde{C_{i}}h}{L}\right\vert
=0.0041, $ and we plot $0.0041\left\vert e^{\left( x_{i+1}-x_{0}\right)
L}-1\right\vert $ as the upper bound (labelled U) on $\left\vert \Delta
_{i}\right\vert .$

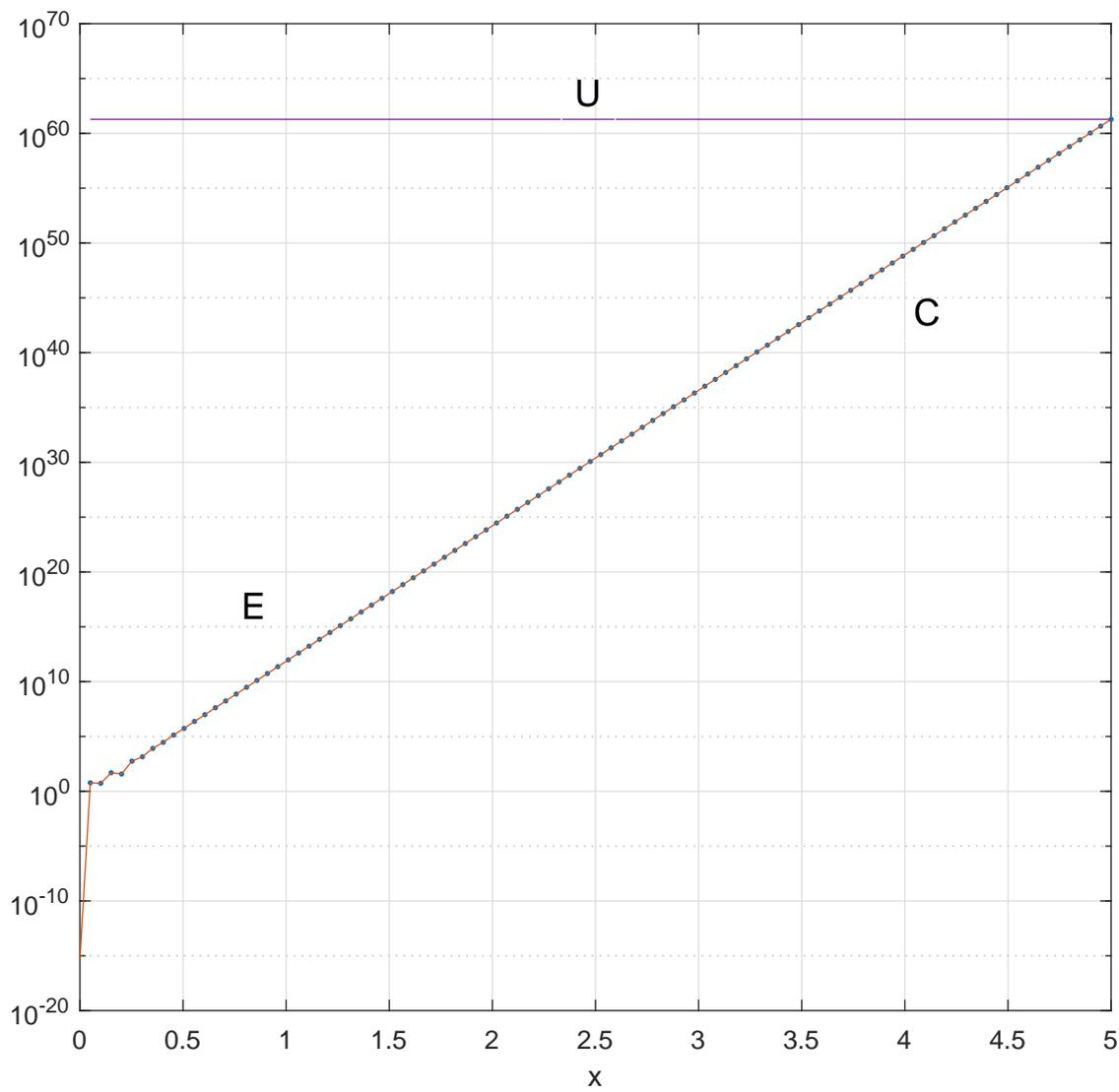
\item \textbf{Figure 2}. We solve (\ref{test eqn}) with $\lambda
=-100,\gamma =-200$ and $h=5\times 10^{-2}$ using the explicit method. The
stepsize is \textit{not} small enough to ensure a stable solution. The
labelling follows that of Figure 1. We estimate $\max \left\vert \frac{%
\widetilde{C_{i}}h}{L}\right\vert =1.9\times 10^{61}.$ As before, $L\ll
0\Rightarrow \left\vert e^{\left( x_{i+1}-x_{0}\right) L}-1\right\vert
\approx 1,$ so that the curve for $\left\vert \frac{\widetilde{C_{i}}h}{L}%
\right\vert $ is superimposed on the curve for $\left\vert \Delta
_{i}\right\vert .$

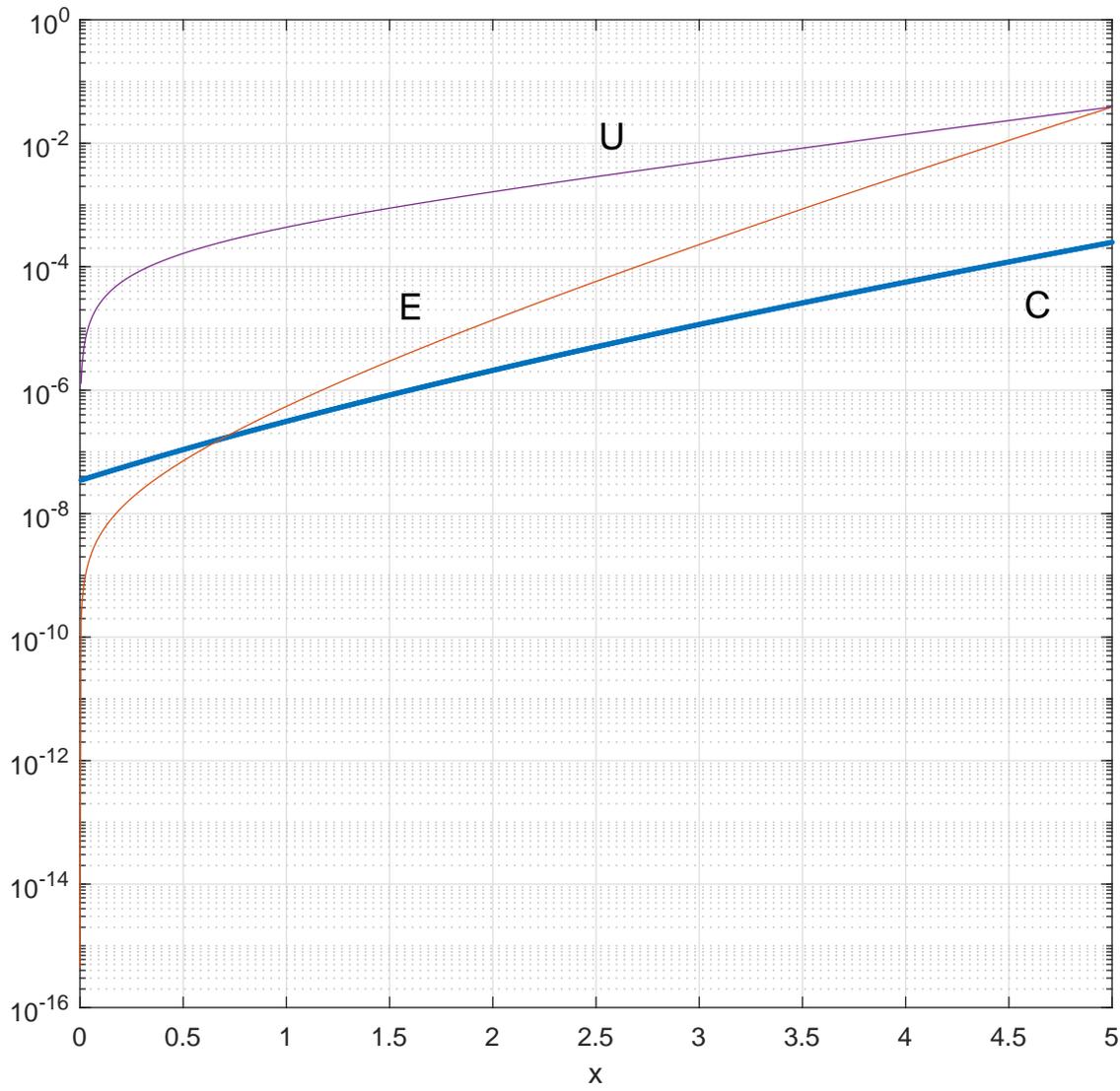
\item \textbf{Figure 3}. We use the explicit method with $\lambda =1,\gamma
=2$ and $h=5\times 10^{-3}.$ We do not have $\left\vert e^{\left(
x_{i+1}-x_{0}\right) L}-1\right\vert \approx 1,$ and so curve C is different
to curve E. We estimate $\max \left\vert \frac{\widetilde{C_{i}}h}{L}%
\right\vert =2.5\times 10^{-4},$ yielding the upper bound U.

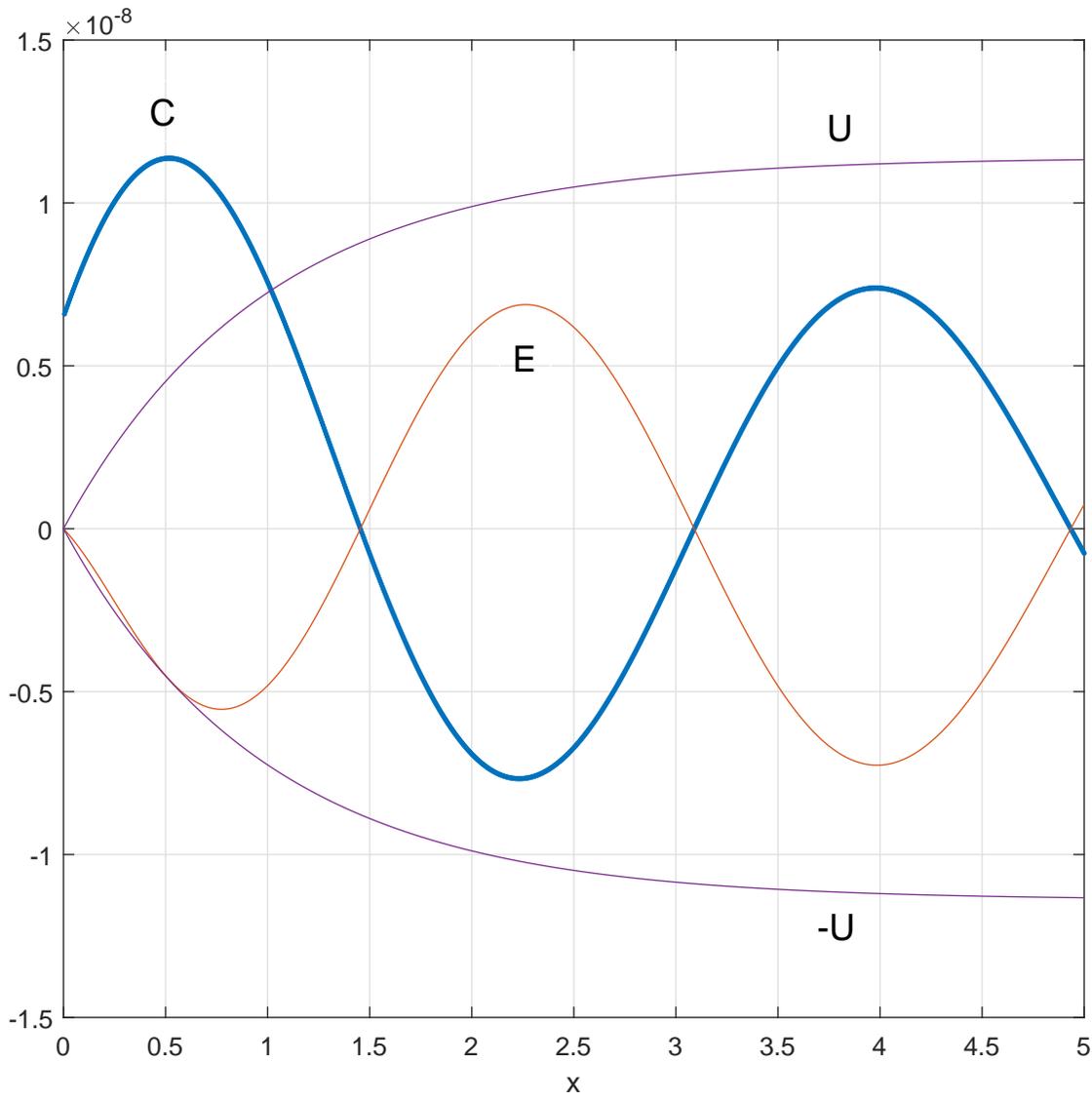
\item \textbf{Figure 4}. Here, we solve the test equation with $\lambda
=-1,\gamma =-2$ and $h=5\times 10^{-3}$ using the implicit method. We show
the signed global error $\Delta _{i},$ and $\frac{\widetilde{C_{i}}h}{L}$.
We see that $\Delta _{i}=0$ when $\frac{\widetilde{C_{i}}h}{L}=0,$ as we
would expect. We estimate $\max \left\vert \frac{\widetilde{C_{i}}h}{L}%
\right\vert =1.14\times 10^{-8},$ yielding the upper and lower bounds (U and
-U). The oscillatory character of the error is due to the oscillatory nature
of the solution.

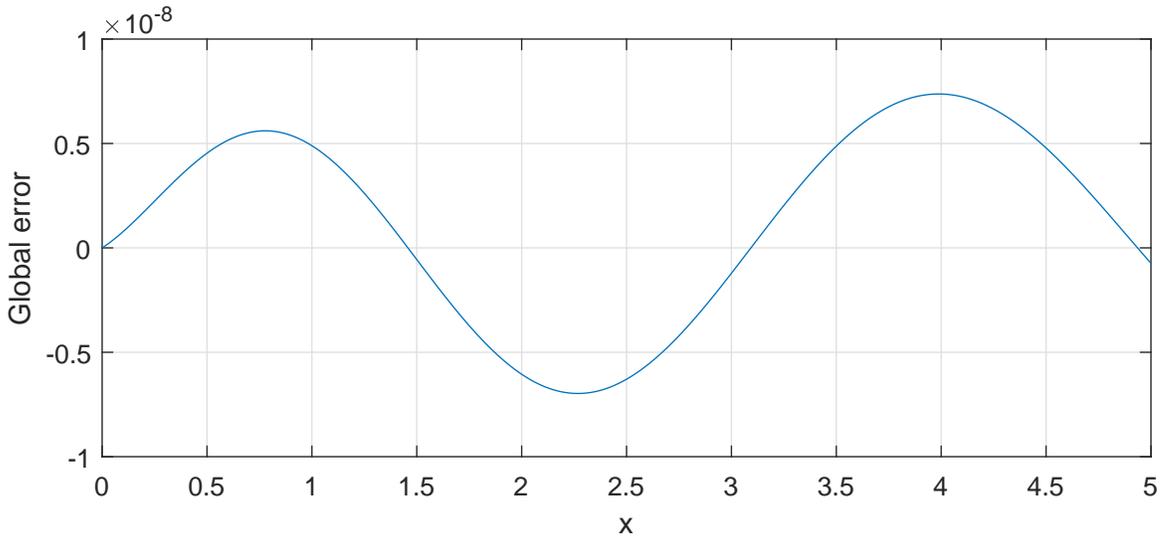
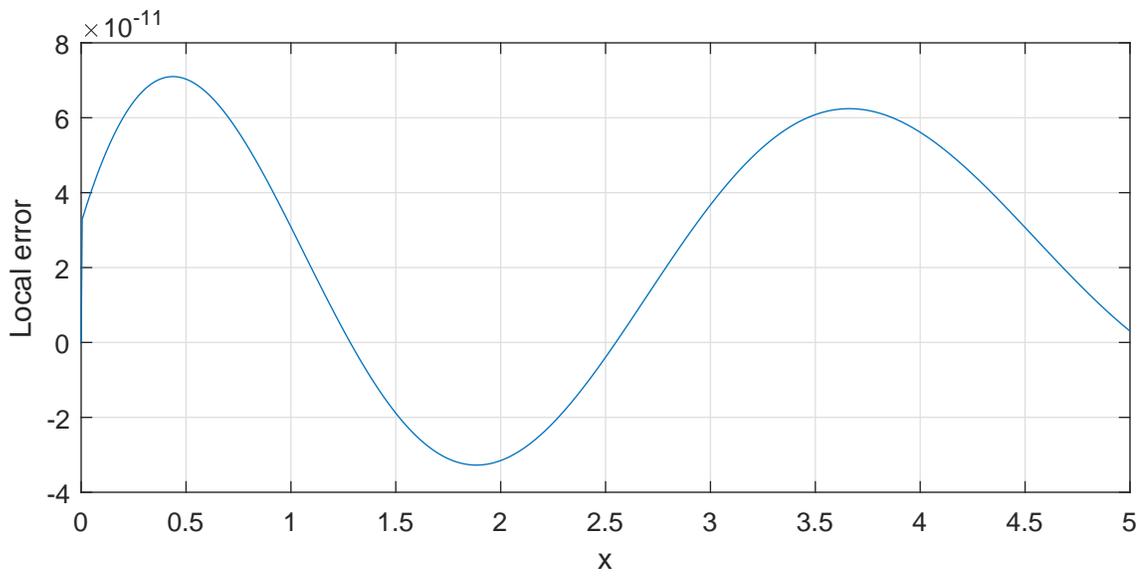
\item \textbf{Figure 5}. We solve the test equation with $\lambda =-1,\gamma
=-2$ and $h=5\times 10^{-3}$ using the explicit method. The upper plot shows
the global error $\Delta _{i},$ and the lower plot shows the local error $%
\varepsilon _{i},$ determined using (\ref{sequence}).
\end{enumerate}

\section{Conclusion}

We have investigated error propagation in an explicit and implicit method
for solving integro-differential equations of the Volterra type. We have
derived upper bounds for the global error, and shown that the global order
for both methods is expected to be first-order. With respect to (\ref%
{problem}), we have considered the case $n=1$. For $n>1,$ we would need to
solve a system of IDEs, and future work would center around error
propagation in such systems - and in systems of IDEs, in general.

\medskip

\section{Appendix}

\subsection{Local order}

The implicit method\relsize{-1}%
\begin{equation*}
y_{i+1}=y_{i}+hf\left( x_{i+1},y_{i+1}\right) +\frac{h^{2}}{2}\left(
\sum\limits_{j=0}^{i+1}2K\left( x_{i+1},y_{j},x_{j}\right) -K\left(
x_{i+1},y_{0},x_{0}\right) -K\left( x_{i+1},y_{i+1},x_{i+1}\right) \right)
\end{equation*}%
\relsize{1}is derived from\relsize{-1}%
\begin{equation*}
\frac{y_{i+1}-y_{i}}{h}=f\left( x_{i+1},y_{i+1}\right) +\frac{h}{2}\left(
\sum\limits_{j=0}^{i+1}2K\left( x_{i+1},y_{j},x_{j}\right) -K\left(
x_{i+1},y_{0},x_{0}\right) -K\left( x_{i+1},y_{i+1},x_{i+1}\right) \right)
\end{equation*}%
\relsize{1}where the LHS is the Euler approximation to the first derivative,
and the second term on the RHS is the composite Trapezium Rule, which models
the integral in (\ref{problem}). As is well-known, the approximation error
in the Euler approximation is $O\left( h\right) $ and in the composite
Trapezium Rule it is $O\left( h^{2}\right) .$ In other words, we have%
\relsize{-1}%
\begin{equation*}
\frac{y_{i+1}-y_{i}}{h}+O\left( h\right) =f\left( x_{i+1},y_{i+1}\right) +%
\frac{h}{2}\left( \sum\limits_{j=0}^{i+1}2K\left( x_{i+1},y_{j},x_{j}\right)
-K\left( x_{i+1},y_{0},x_{0}\right) -K\left( x_{i+1},y_{i+1},x_{i+1}\right)
\right) +O\left( h^{2}\right)
\end{equation*}%
\relsize{1}which gives\relsize{-1}%
\begin{equation*}
y_{i+1}-y_{i}+O\left( h^{2}\right) =hf\left( x_{i+1},y_{i+1}\right) +\frac{%
h^{2}}{2}\left( \sum\limits_{j=0}^{i+1}2K\left( x_{i+1},y_{j},x_{j}\right)
-K\left( x_{i+1},y_{0},x_{0}\right) -K\left( x_{i+1},y_{i+1},x_{i+1}\right)
\right) +O\left( h^{3}\right) .
\end{equation*}%
\relsize{1}The sum of the $O\left( h^{2}\right) $ term and the $O\left(
h^{3}\right) $ is the local error $\varepsilon _{i+1},$ so that 
\begin{equation*}
\varepsilon _{i+1}=\underset{\varepsilon _{i+1}^{D}}{\underbrace{O\left(
h^{2}\right) }}+\underset{\varepsilon _{i+1}^{Q}}{\underbrace{O\left(
h^{3}\right) }}=O\left( h^{2}\right)
\end{equation*}%
for the implicit method.

For the explicit method\relsize{-1} 
\begin{equation*}
y_{i+1}=y_{i}+hf\left( x_{i},y_{i}\right) +\frac{h^{2}}{2}\left(
\sum\limits_{j=0}^{i}2K\left( x_{i},y_{j},x_{j}\right) -K\left(
x_{i},y_{0},x_{0}\right) -K\left( x_{i},y_{i},x_{i}\right) \right)
\end{equation*}%
\relsize{1}the composite Trapezium Rule is written\relsize{-1}%
\begin{align*}
& \frac{h}{2}\left( \sum\limits_{j=0}^{i+1}2K\left(
x_{i+1},y_{j},x_{j}\right) -K\left( x_{i+1},y_{0},x_{0}\right) -K\left(
x_{i+1},y_{i+1},x_{i+1}\right) \right) +O\left( h^{2}\right) \\
& =\frac{h}{2}\left( \sum\limits_{j=0}^{i}2K\left( x_{i},y_{j},x_{j}\right)
-K\left( x_{i},y_{0},x_{0}\right) -K\left( x_{i},y_{i},x_{i}\right) \right) +%
\frac{h}{2}\left( K\left( x_{i},y_{i},x_{i}\right) -K\left(
x_{i+1},y_{i+1},x_{i+1}\right) \right) +O\left( h^{2}\right) \\
& =\frac{h}{2}\left( \sum\limits_{j=0}^{i}2K\left( x_{i},y_{j},x_{j}\right)
-K\left( x_{i},y_{0},x_{0}\right) -K\left( x_{i},y_{i},x_{i}\right) \right) +%
\underset{O\left( h\right) }{\underbrace{O\left( h\right) +O\left(
h^{2}\right) }}.
\end{align*}%
\relsize{1}Consequently, the local error for the explicit method has the
form 
\begin{equation*}
\varepsilon _{i+1}=\underset{\varepsilon _{i+1}^{D}}{\underbrace{O\left(
h^{2}\right) }}+\underset{\varepsilon _{i+1}^{Q}}{\underbrace{O\left(
h^{2}\right) }}=O\left( h^{2}\right) .
\end{equation*}

It is worth noting that, for both methods, $\varepsilon _{i+1}\rightarrow 0$
as $h\rightarrow 0,$ implying \textit{consistency}.

\subsection{Roundoff}

For completeness' sake, we can include a roundoff term $\mu _{i}$ in each
local error, as in%
\begin{align*}
\widetilde{\varepsilon }_{i}& \rightarrow \widetilde{\varepsilon }_{i}+\mu
_{i} \\
\widetilde{\widetilde{\varepsilon }}_{i}& \rightarrow \widetilde{\widetilde{%
\varepsilon }}_{i}+\mu _{i}.
\end{align*}%
This will simply lead to terms of the form 
\begin{align*}
& \left( \frac{\widetilde{C}h}{L}+\frac{\mu _{i}}{hL}\right) \left(
e^{\left( x_{i+1}-x_{0}\right) L}-1\right) \\
& \left( \frac{\widetilde{\widetilde{C}}h}{L}+\frac{\mu _{i}}{hL}\right)
\left( e^{\left( x_{i+1}-x_{0}\right) L}-1\right)
\end{align*}%
in the upper bounds derived earlier. If $\left\vert hL\right\vert \ll 1,$
the roundoff component could be significant, but this is unlikely -
particularly in the context of modern computing, where numerical precision
can be controlled via the \textit{variable precision arithmetic} (VPA)
capabilities of current software. In other words, we can effectively make $%
\mu _{i}$ as small as we need it to be. This may occur at the cost of
computational efficiency, but that is the way of things.

\end{document}